\documentstyle[12pt]{article}

\font\twlgot =eufm10 scaled \magstep1 \font\egtgot =eufm8
\font\sevgot =eufm7 \font\twlmsb =msbm10 scaled \magstep1
\font\egtmsb =msbm8 \font\sevmsb =msbm7

\newfam\gotfam
\def\pgot{\fam\gotfam\twlgot}
\textfont\gotfam\twlgot \scriptfont\gotfam\egtgot
\scriptscriptfont\gotfam\sevgot
\def\got{\protect\pgot}
\newfam\msbfam
\textfont\msbfam\twlmsb \scriptfont\msbfam\egtmsb
\scriptscriptfont\msbfam\sevmsb
\def\Bbb{\protect\pBbb}
\def\pBbb{\relax\ifmmode\expandafter\Bb\else\typeout{You cann't use
Bbb in text mode}\fi}
\def\Bb #1{{\fam\msbfam\relax#1}}

\newcommand{\gQ}{{\got T}}

\newcommand{\gA}{{\got A}}

\newcommand{\gd}{{\got d}}

\def\thebibliography#1{\section*{References}\list
  {[\arabic{enumi}]}{\settowidth\labelwidth{#1}\leftmargin\labelwidth
    \advance\leftmargin\labelsep
    \usecounter{enumi}}
    \def\newblock{\hskip .11em plus .33em minus .07em}
    \sloppy\clubpenalty4000\widowpenalty4000
    \sfcode`\.=1000\relax}

\def\op#1{\mathop{\fam0 #1}\limits}

\newcommand{\Ker}{{\rm Ker\,}}
\newcommand{\im}{{\rm Im\,}}

\newcommand{\beq}{\begin{equation}}
\newcommand{\eeq}{\end{equation}}
\newcommand{\ben}{\begin{eqnarray}}
\newcommand{\een}{\end{eqnarray}}
\newcommand{\be}{\begin{eqnarray*}}
\newcommand{\ee}{\end{eqnarray*}}
\newcommand{\bea}{\begin{eqalph}}
\newcommand{\eea}{\end{eqalph}}
\newcommand{\cA}{{\cal A}}

\newcommand{\cP}{{\cal P}}

\newcommand{\cE}{{\cal E}}

\newcommand{\cS}{{\cal S}}
\newcommand{\cC}{{\cal C}}
\newcommand{\cO}{{\cal O}}

\newcommand{\dl}{\delta}
\newcommand{\la}{\lambda}
\newcommand{\La}{\Lambda}
\newcommand{\f}{\phi}

\newcommand{\m}{\mu}

\newcommand{\G}{\Gamma}

\newcommand{\vt}{\vartheta}

\newcommand{\up}{\upsilon}

\newcommand{\si}{\sigma}
\newcommand{\Si}{\Sigma}
\newcommand{\w}{\wedge}

\newcommand{\ol}{\overline}
\newcommand{\dr}{\partial}
\newcommand{\ar}{\op\longrightarrow}
\newcommand{\llr}{\op\longleftarrow}
\newcommand{\lto}{\leftarrow}

\newcommand{\rdr}{\stackrel{\leftarrow}{\dr}{}}

\newcounter{eqalph}
\newcounter{equationa}
\newcounter{remark}
\newcounter{example}
\newcounter{theorem}
\newcounter{proposition}
\newcounter{lemma}
\newcounter{corollary}
\newcounter{definition}
\newenvironment{eqalph}{\stepcounter{equation}
\setcounter{equationa}{\value{equation}} \setcounter{equation}{0}

\begin{eqnarray}}{\end{eqnarray}\setcounter{equation}{\value{equationa}}}

\def\theremark{\arabic{remark}}

\def\thetheorem{\arabic{theorem}}

\newenvironment{proof}{
{\bf Proof.}}{\hfill $\Box$ \medskip}
\newenvironment{rem}{\refstepcounter{remark}\medskip\noindent{\bf
Remark \theremark.}}{\medskip}

\newenvironment{theo}{\refstepcounter{theorem}\medskip\noindent
{\bf Theorem \thetheorem.}\it }{\medskip}
\newenvironment{prop}{\refstepcounter{theorem}\medskip\noindent
{\bf Proposition \thetheorem.}\it }{\medskip}
\newenvironment{lem}{\refstepcounter{theorem}\medskip\noindent
{\bf Lemma \thetheorem.}\it }{\medskip}

\newenvironment{defi}{\refstepcounter{theorem}\medskip\noindent
{\bf Definition \thetheorem.} \it }{\medskip}

\newcommand{\mar}[1]{}

\begin{document}
\hbox{}

\begin{center}


{\large \bf Noether identities of a generic differential operator.
The Koszul--Tate complex}
\bigskip

{\sc G. SARDANASHVILY}
\medskip

\begin{small}

 {\it Department of Theoretical Physics,
Moscow State University}



\end{small}
\end{center}

\bigskip

\begin{small}
\noindent
 Given a generic Lagrangian system, its Euler--Lagrange
operator obeys Noether identities which need not be independent,
but satisfy first-stage Noether identities, and so on. This
construction is generalized to arbitrary differential operators on
a smooth fiber bundle. Namely, if a certain necessary and
sufficient condition holds, one can associate to a differential
operator the exact chain complex with the boundary operator whose
nilpotency restarts all the Noether identities characterizing the
degeneracy of an original differential operator.


\end{small}

\section{Introduction}

Given a generic Lagrangian system on a fiber bundle, its
Euler--Lagrange operator $\cE$ obeys the Noether identities. They
need not be independent, but satisfy the first-stage Noether
identities, which in turn are subject to the second-order Noether
ones, and so on. The hierarchy of these Noether identities
characterizes the degeneracy of a Lagrangian system in full.
Noether's second theorem states the relation between Noether
identities and gauge symmetries of a Lagrangian
\cite{jpa05,jmp05}, but Noether identities can be introduced
without regarding $\cE$ as an Euler--Lagrange operator. Therefore,
one can extend the notion of Noether identities to a generic
differential operator on a fiber bundle. Our goal is the
following.

Let $Y\to X$ be a smooth fiber bundle and $J^rY$, $r=1,\ldots,$
the jet manifolds of its sections (the index $r=0$ further stands
for $Y$). Let $E\to X$ be a vector bundle. A $E$-valued $r$-order
differential operator $\cE$ on $Y$ is conventionally defined as a
bundle morphism of $J^rY$ to $E$ over $X$ \cite{kras,bry}. We
associate to $\cE$ the chain complex  (\ref{v42}) whose boundaries
vanish on $\Ker\cE$ (Proposition \ref{v203}). It is a complex of a
certain ring $\cS_\infty[E;Y]$ of Grassmann-graded functions and
their jets on the infinite order jet manifold $J^\infty Y$ of $Y$.
For our purpose, this complex can be replaced with the short
zero-exact complex $\cS_\infty[E;Y]_{\leq 2}$ (\ref{v042}). Recall
that a chain complex is called $r$-exact if its homology of degree
$k\leq r$ is trivial.

Noether identities of a differential operator $\cE$ are defined as
nontrivial elements of the first homology $H_1(\dl)$ of the
complex (\ref{v042}) (Definition \ref{v121}). If this homology is
finitely generated, the complex (\ref{v042}) can be extended to
the one exact-complex $\cP_\infty[E;Y;E_0]_{\leq 3}$ (\ref{v66})
with the boundary operator $\dl_0$ (\ref{v204}) whose nilpotency
conditions are equivalent to the above-mentioned Noether
identities (Proposition \ref{v137}). First-stage Noether
identities are defined as two-cycles of this complex. They are
trivial if two-cycles are boundaries, but the converse need not be
true. Trivial first-stage Noether identities are boundaries iff a
certain homology condition (called the two-homology regularity
condition) holds (Proposition \ref{v134}). In this case, the
first-stage Noether identities are identified to nontrivial
elements of the second homology of the complex (\ref{v66}). If
this homology is finitely generated, the complex (\ref{v66}) is
extended to the two-exact complex $\cP_\infty[E_1E;Y;E_0]_{\leq
4}$ (\ref{v87}) with the boundary operator $\dl_1$ (\ref{v205})
whose nilpotency conditions are equivalent to the Noether and
first-stage Noether identities (Proposition \ref{v139}). If the
third homology of this complex is not trivial, the second-stage
Noether identities exist, and so on. Iterating the arguments, we
come to the following.

Let we have the $(N+1)$-exact complex $\cP_\infty\{N\}_{\leq N+3}$
(\ref{v94}) such that: (i) the nilpotency conditions of its
boundary operator $\dl_N$ (\ref{v92}) reproduce Noether and
$k$-stage Noether identities for $k\leq N$, (ii) the
$(N+1)$-homology regularity condition (Definition \ref{v155})
holds. Then the $(N+1)$-stage Noether identities are defined as
$(N+2)$-cycles of this complex. They are trivial if cycles are
boundaries, while the converse is true iff the $(N+2)$-homology
regularity condition is satisfied. In this case, $(N+1)$-stage
Noether identities are identified to nontrivial elements of the
$(N+2)$-homology of the complex (\ref{v94}) (item (i) of Theorem
\ref{v163}). If this homology is finitely generated, this complex
is extended to the $(N+2)$-exact complex $\cP_\infty\{N+1\}_{\leq
N+4}$ (\ref{v171}) with the boundary operator $\dl_{N+1}$
(\ref{v170}) whose nilpotency restarts all the Noether identities
up to stage $(N+1)$  (item (ii) of Theorem \ref{v163}).

This iteration procedure results in the exact Koszul--Tate complex
with the boundary operator whose nilpotency conditions reproduce
all Noether and higher Noether identities characterizing the
degeneracy of a differential operator $\cE$

\section{The ring of Grassmann-graded functions and their jets}

All chain complexes considered in the article are complexes of
certain rings of Grassmann-graded function and their jets on the
infinite order jet manifold $J^\infty Y$ of a fiber bundle $Y\to
X$. Let us describe such a ring.

Recall that $J^\infty Y$, is the projective limit
$(\pi^\infty_r:J^\infty Y\to J^rY)$ of the inverse system of jet
manifolds
\mar{5.10}\beq
X\op\longleftarrow^\pi Y\op\longleftarrow^{\pi^1_0} J^1Y
\longleftarrow \cdots J^{r-1}Y \op\longleftarrow^{\pi^r_{r-1}}
J^rY\longleftarrow\cdots, \label{5.10}
\eeq
of $J^\infty Y$, where $\pi^r_{r-1}$ are affine bundles. It is a
Fr\'echet manifold.  A bundle atlas $\{(U_Y;x^\la,y^i)\}$ of $Y\to
X$ induces the coordinate atlas
\mar{jet1}\beq
\{((\pi^\infty_0)^{-1}(U_Y); x^\la, y^i_\La)\}, \qquad
{y'}^i_{\la+\La}=\frac{\dr x^\m}{\dr x'^\la}d_\m y'^i_\La, \qquad
0\leq|\La|, \label{jet1}
\eeq
where $\La=(\la_1...\la_k)$, $\la+\La=(\la\la_1...\la_k)$ are
symmetric multi-indices and
\be
d_\la = \dr_\la + \op\sum_{0\leq|\La|} y^i_{\la+\La}\dr_i^\La,
\qquad d_\La=d_{\la_1}\circ\cdots\circ d_{\la_k},
\ee
are total derivatives. We further assume that the cover
$\{\pi(U_Y)\}$ of $X$ is also the cover of atlases of all vector
bundles over $X$ in question. The inverse system (\ref{5.10})
yields the direct system
\be
C^\infty(X)\op\longrightarrow^{\pi^*} C^\infty(Y)
\op\longrightarrow^{\pi^1_0{}^*} C^\infty(J^1Y) \ar\cdots
C^\infty(J^{r-1}Y) \op\longrightarrow^{\pi^r_{r-1}{}^*}
 \cdots
\ee
of rings of smooth real functions on jet manifolds $J^rY$ with
respect to the pull-back monomorphisms $\pi^r_{r-1}{}^*$. Its
direct limit is the ring $\cO_\infty Y$ of all smooth real
functions on finite order jet manifolds modulo the pull-back
identification.

Let us extend the ring $\cO_\infty Y$ to a ring of graded
functions on graded manifolds whose bodies are jet manifolds
$J^rY$ of $Y$ \cite{jmp05,cmp04}. We restrict our consideration to
graded manifolds $(Z,\gA)$ with structure sheaves $\gA$ of
Grassmann algebras of finite rank, and refer to the following
Serre--Swan theorem.

\begin{theo} \label{v0}
Let $Z$ be a smooth manifold. A Grassmann algebra $\cA$ over the
ring $C^\infty(Z)$ of smooth real functions on $Z$ is isomorphic
to the Grassmann algebra of graded functions on a graded manifold
with a body $Z$ iff it is the exterior algebra of some projective
$C^\infty(Z)$-module of finite rank.
\end{theo}

\begin{proof} The proof follows at once from the Batchelor theorem
\cite{bart} and the Serre-Swan theorem generalized to an arbitrary
smooth manifold \cite{book05,ren}. The Batchelor theorem states
that any graded manifold $(Z,\gA)$ with a body $Z$ is isomorphic
to the one $(Z,\gA_Q)$ with the structure sheaf $\gA_Q$ of germs
of sections of the exterior bundle
\be
\w Q^*=\Bbb R\op\oplus_Z Q^*\op\oplus_Z\op\w^2
Q^*\op\oplus_Z\cdots,
\ee
where $Q^*$ is the dual of some vector bundle $Q\to Z$. We agree
to call $(Z,\gA_Q)$ the simple graded manifold with the structure
vector bundle $Q$. Its structure ring $\cA_Q$ of graded functions
(sections of $\gA_Q$) is the $\Bbb Z_2$-graded exterior algebra of
the $C^\infty(Z)$-module of sections of $Q^*\to Z$. By virtue of
the Serre--Swan theorem, a $C^\infty(Z)$-module is isomorphic to
the module of sections of a smooth vector bundle over $Z$ iff it
is a projective module of finite rank.
\end{proof}

\begin{rem} With respect to bundle coordinates
$(z^A,q^a)$ on $Q$ and the corresponding fiber basis $\{c^a\}$ for
$Q^*\to X$, graded functions read
\mar{v23}\beq
f=\op\sum_{k=0} \frac1{k!}f_{a_1\ldots a_k}c^{a_1}\cdots c^{a_k},
\qquad f\in C^\infty(Z), \label{v23}
\eeq
where we omit the symbol of the exterior product of elements
$c^a$. Let $u\in \gd\cA_Q$ be a graded derivation of the $\Bbb
R$-ring $\cA_Q$. Due to the canonical splitting $VQ= Q\times Q$,
the fiber basis $\{\dr_a\}$ for vertical tangent bundle $VQ\to Q$
of $Q\to Z$ is the dual of $\{c^a\}$ and, therefore, $u$ takes the
local form $u= u^A\dr_A + u^a\dr_a$, where $u^A, u^a$ are local
graded functions \cite{bart}. It acts on graded functions
(\ref{v23}) by the rule
\be
 u(f_{a\ldots b}c^a\cdots c^b)=u^A\dr_A(f_{a\ldots
b})c^a\cdots c^b +u^d f_{a\ldots b}\dr_d\rfloor (c^a\cdots c^b).
\ee
\end{rem}

Given a vector bundle $E\to X$, let us consider the simple graded
manifold $(J^rY,\gA_{E_r})$ whose body is $J^rY$ and the structure
bundle is the pull-back
\be
E_r=J^rY\op\times_X J^rE
\ee
onto $J^rY$ of the jet bundle $J^rE\to X$, which is a vector
bundle. There is an epimorphism of graded manifolds
\be
(J^{r+1}Y,\gA_{E_{r+1}}) \to (J^rY,\gA_{E_r}).
\ee
It consists of the open surjection $\pi^{r+1}_r$ and the sheaf
monomorphism $\pi_r^{r+1*}:\gA_{E_r}\to \gA_{E_{r+1}}$, where
$\pi_r^{r+1*}\gA_{E_r}$ is the pull-back onto $J^{r+1}Y$ of the
topological fiber bundle $\gA_{E_r}\to J^rY$. These sheaf
monomorphisms yield monomorphisms of the corresponding canonical
presheaves $\ol\gA_{E_r}\to \ol\gA_{E_{r+1}}$, which make up a
direct system
\mar{v15}\beq
\ol\gA_{Y\times E}\ar \ol\gA_{E_1} \ar\cdots \ol\gA_{E_r}
\ar\cdots, \label{v15}
\eeq
and the monomorphisms of graded commutative rings $\cA_{E_r} \to
\cA_{E_{r+1}}$ assembled into the direct system
\mar{v6}\beq
\cA_{Y\times E}\ar \cA_{E_1}\ar\cdots \cA_{E_r}\ar\cdots.
\label{v6}
\eeq
A direct limit of this direct system is a graded commutative ring
$\cS_\infty[E;Y]$ of all graded functions $f\in \cA_{E_r}$ on jet
manifolds $J^rY$ modulo monomorphisms $\pi_r^{r+1*}$. The
monomorphisms $C^\infty(J^rY)\to \cA_{E_r}$ provide the
monomorphism $\cO_\infty Y\to \cS_\infty[E;Y]$, while the body
epimorphisms $\cA_{E_r}\to C^\infty(J^rY)$ yield the epimorphism
$\cS_\infty[E;Y]\to \cO_\infty Y$.

One can think of elements of $\cS_\infty[E;Y]$ as being graded
functions on $J^\infty Y$ as follows. A direct limit of the direct
system of presheaves (\ref{v15}) is a presheaf on the infinite
order jet manifold $J^\infty Y$. Let $\gQ_\infty[E;Y]$ be the
sheaf of germs of this presheaf. The structure module
$\G(\gQ_\infty[E;Y])$ of sections of $\gQ_\infty[E;Y]$ is a ring
such that, given an element $f\in \G(\gQ_\infty[E;Y])$ and a point
$z\in J^\infty Y$, there exist an open neighbourhood $U$ of $z$
and a graded function $f^{(k)}$ on some finite order jet manifold
$J^kY$ so that $f|_U= \pi^{\infty*}_k f^{(k)}|_U$. In particular,
there is the monomorphism $\cS_\infty[E;Y]
\to\G(\gQ_\infty[E;Y])$. Due to this monomorphism, one can
restrict $\cS_\infty[E;Y]$ to the coordinate chart (\ref{jet1}),
and say that $\cS_\infty[E;Y]$ as an $\cO_\infty Y$-ring is
locally generated by the elements $c^a_\La$, where $\{c^a\}$ is a
a local fiber basis for $E^*$ over $\pi_\infty(U)$. We agree to
call $(y^i_\La,c^a_\La)$, $0\leq |\La|$, the local basis for
$\cS_\infty[E;Y]$. Let the collective symbol $s^A_\La$ stand for
its elements.

\begin{rem}
One can think of $c^a_\La$ as being jets of graded functions $c^a$
\cite{jmp05,cmp04}. This definition differs from the notion of
jets of a graded commutative ring \cite{book05} and that of jets
of a graded fiber bundle \cite{hern}, but reproduces the heuristic
notion of jets of odd fields in Lagrangian field theory
\cite{barn}.
\end{rem}

\begin{rem} \label{v201} \mar{v201} Let $V$, $V'$ and $V''$ be vector
bundles over $X$.
Let us consider the ring
\be
\cS_\infty[V\op\times_X E;Y\op\times_X V'\op\times_X V''],
\ee
and its subring of graded functions which are polynomial in fiber
coordinates of the vector bundle $Y\op\times_X V'\op\times_X
V''\to Y$. We denote the latter by $\cP_\infty[VE;Y;V'V'']$. One
can think of its elements as being graded functions on $J^\infty
Y$, too.
\end{rem}

Let $\vt\in\gd \cS_\infty[E;Y]$ be a graded derivation of the
$\Bbb R$-ring $\cS_\infty[E;Y]$ \cite{jmp05,cmp04}. With respect
to the local basis $(s^A_\La)$, it takes the form
\be
 \vt=\vt^\la d_\la + (\vt^A\dr_A
+\op\sum_{|\La|>0} d_\La\vt^A\dr_A^\La),
\ee
where the tuple of graded derivations $\{\dr^\La_a\}$ is the dual
of the tuple $\{c^a_\La\}$ of generating elements of the
$\cO_\infty[F;Y]$-algebra $\cS_\infty[E;Y]$, and $\vt^\la$,
$\vt^A$ are local graded functions. We further restrict our
consideration to vertical contact graded derivations
\mar{0672}\beq
\vt=\op\sum_{0\leq|\La|} d_\La\up^A\dr_A^\La. \label{0672}
\eeq
Such a derivation is completely determined by its first summand
\mar{0673}\beq
\up=\up^A(x^\la,s^A_\La)\dr_A, \qquad 0\leq|\La|\leq k,
\label{0673}
\eeq
called a generalized graded vector field. For the sake of
simplicity, the common symbol $\up$ further stands both for a
contact graded derivation (\ref{0672}) and a generalized vector
field (\ref{0673}), which is also called a graded derivation. A
graded derivation $\up$ is said to be nilpotent if $\up(\up(f))=0$
for any graded function $f\in \cS_\infty[E;Y]$. One can show that
$\up$ is nilpotent only if it is odd and iff all $\up^A$ obey the
equality
\mar{0688}\beq
\up(\up^A)=\op\sum_{0\leq|\Si|} \up^B_\Si\dr^\Si_B(\up^A)=0.
\label{0688}
\eeq

\section{Noether identities}

Let $\cE$ be a $E$-valued differential operator on a smooth fiber
bundle $Y\to X$.

\begin{prop} \label{v203} \mar{v203}
One can associate to $\cE$ a chain complex whose boundaries vanish
on $\Ker\cE$.
\end{prop}

\begin{proof}
A differential operator $\cE$ is locally represented by a set of
functions $\cE^a\in \cO_\infty Y$, possessing the corresponding
coordinate transformation law. Let us provide the ring
$\cS_\infty[E;Y]$ with the nilpotent graded derivation
\mar{vv}\beq
\dl=\rdr_a\cE^a, \label{vv}
\eeq
whose definition is independent of the choice of a local basis. It
is convenient to deal with a graded derivation $\dl$ (\ref{vv})
acting on graded functions on the right by the rule
\be
\dl(f\w f')=(-1)^{[f']}\dl(f)\w f'+ f\w\dl(f'),
\ee
where $[f']$ denotes the Grassmann parity. We call $\dl$
(\ref{vv}) the Koszul--Tate differential. With this differential,
the ring $\cS_\infty[E;Y]$ is split into the chain complex
\mar{v42}\beq
0\lto \cO_\infty Y \llr^\dl \cS_\infty[E;Y]_1\cdots \llr^\dl
\cS_\infty[E;Y]_k \cdots \label{v42}
\eeq
graded in polynomials of the odd elements $c^a_\La$. Following the
physical literature \cite{jmp05,barn}, we assign to $c^a_\La$ the
antifield number Ant$[c^a_\La]=1$. It is readily observed that the
boundaries of the complex (\ref{v42}) vanish on $\Ker \cE$.
\end{proof}

Note that homology groups $H_*(\dl)$ of the complex (\ref{v42})
are $\cO_\infty Y$-modules, but these modules fail to be
torsion-free. Indeed, given a cycle $\Phi\in \cS_\infty[E;Y]_k$
and an element $f=\dl\si$ of the ring $\cO_\infty Y\subset
\cS_\infty[E;Y]$, we obtain that $f\Phi=\dl(\si \Phi)$ is a
boundary. Therefore, one can not apply the K\"unneth formula to
the homology of this complex, though any its term
$\cS_\infty[E;Y]_k$ is isomorphic to the graded commutative
$k$-tensor product of the $\cO_\infty Y$-module
$\cS_\infty[E;Y]_1$.

The homology $H_0(\dl)$ of the complex (\ref{v42}) is not trivial,
but this homology and the higher ones $H_{k\geq 2}(\dl)$ are not
essential for our consideration. Therefore, we replace the complex
(\ref{v42}) with the finite one
\mar{v042}\beq
0\lto \im\dl \llr^\dl \cS_\infty[E;Y]_1 \llr^\dl \cS_\infty[E;Y]_2
\label{v042}
\eeq
of graded functions of antifield number $k\leq 2$. It is exact at
$\im\dl$, and  its first homology coincides with that of  the
complex (\ref{v42}). Let us consider this homology.

A generic one-chain of the complex (\ref{v042}) takes the form
\mar{0712}\beq
\Phi= \op\sum_{0\leq|\La|} \Phi_a^\La c^a_\La, \qquad
\Phi_a^\La\in \cO_\infty Y, \label{0712}
\eeq
and the cycle condition $\dl \Phi=0$ reads
\mar{0713}\beq
\op\sum_{0\leq|\La|} \Phi_a^\La d_\La \cE^a =0. \label{0713}
\eeq
One can think of this equality as being a reduction condition on a
differential operator $\cE$. Conversely, any reduction condition
of form (\ref{0713}) comes from some cycle (\ref{0712}). The
reduction condition (\ref{0713}) is trivial  if a cycle is a
boundary, i.e., it takes the form
\mar{v44}\beq
\Phi= \op\sum_{0\leq|\La|,|\Si|} T_{ab}^{\La\Si}d_\Si\cE^b
c^a_\La, \qquad T_{ab}^{\La\Si}=- T_{ba}^{\Si\La}. \label{v44}
\eeq
If $\cE$ is an Euler--Lagrange operator of some Lagrangian system
on a fiber bundle $Y\to X$, the nontrivial reduction condition
(\ref{0713}) is a Noether identity \cite{jpa05,jmp05}. Therefore,
we come to the following definition.

\begin{defi} \label{v121} \mar{v121}
A differential operator $\cE$ is called degenerate if the homology
$H_1(\dl)$ of the complex (\ref{v44}) (or (\ref{0712})) is not
trivial. We agree to call a cycle condition (\ref{0713}) the
Noether identity.
\end{defi}

One can say something more if the $\cO_\infty Y$-module $H_1(\dl)$
is finitely generated, i.e., it possesses the following particular
structure. There are elements $\Delta\in H_1(\dl)$ making up a
projective $C^\infty(X)$-module $\cC_{(0)}$ of finite rank which,
by virtue of the Serre--Swan theorem, is isomorphic to the module
of sections of some vector bundle $E_0\to X$. Let $\{\Delta^r\}$
be local bases for this $C^\infty(X)$-module. Then every element
$\Phi\in H_1(\dl)$ factorizes
\mar{v63,71}\ben
&& \Phi= \op\sum_{0\leq|\Xi|} G_r^\Xi d_\Xi \Delta^r, \qquad
G_r^\Xi\in
\cO_\infty Y, \label{v63}\\
&&\Delta^r=\op\sum_{0\leq|\La|} \Delta_a^{\La r} c^a_\La, \qquad
\Delta_a^{\La r}\in \cO_\infty Y, \label{v71}
\een
via elements of $\cC_{(0)}$, i.e., any Noether identity
(\ref{0713}) is a corollary of the Noether identities
\mar{v64}\beq
 \op\sum_{0\leq|\La|} \Delta_a^{\La r} d_\La \cE^a=0.
\label{v64}
\eeq
Clearly, the factorization (\ref{v63}) is independent of
specification of local bases $\{\Delta^r\}$. We say that the
Noether identities  (\ref{v64}) are complete, and call $\Delta\in
\cC_{(0)}$ the Noether operators. Note that, being representatives
of $H_1(\dl)$, the graded functions $\Delta^r$ (\ref{v71}) are not
$\dl$-exact.

\begin{prop} \label{v137} \mar{v137}
If the homology $H_1(\dl)$ of the complex (\ref{v042}) is finitely
generated, this complex  can be extended to a one-exact complex
with a boundary operator whose nilpotency conditions (\ref{0688})
are just complete Noether identities (see the complex (\ref{v66})
below).
\end{prop}

\begin{proof}
Let us consider the ring $\cP_\infty[E;Y;E_0]$ of graded functions
on $J^\infty Y$ (see Remark \ref{v201}). It possesses local bases
$\{y^i_\La, c^a_\La, c^r_\La\}$, where $[c^r_\La]=0$ and
Ant$[c^r_\La]=2$. This ring is provided with the nilpotent graded
derivation
\mar{v204}\beq
\dl_0=\dl + \rdr_r\Delta^r, \label{v204}
\eeq
called the extended Koszul--Tate differential. Its nilpotency
conditions (\ref{0688}) are equivalent to the complete Noether
identities (\ref{v64}). Then the module $\cP_\infty[E;Y;E_0]_{\leq
3}$ of graded functions of antifield number $k\leq 3$ is split
into the chain complex
\mar{v66}\beq
0\lto \im\dl \llr^\dl \cS_\infty[E;Y]_1\llr^{\dl_0}
\cP_\infty[E;Y;E_0]_2\llr^{\dl_0} \cP_\infty[E;Y;E_0]_3.
\label{v66}
\eeq
Let $H_*(\dl_0)$ denote its homology. We have
$H_0(\dl_0)=H_0(\dl)=0$. Furthermore, any one-cycle $\Phi$ up to a
boundary takes the form (\ref{v63}) and, therefore, it is a
$\dl_0$-boundary
\be
\Phi= \op\sum_{0\leq|\Si|} G_r^\Xi d_\Xi \Delta^r
=\dl_0(\op\sum_{0\leq|\Si|} G_r^\Xi c^r_\Xi).
\ee
Hence, $H_1(\dl_0)=0$, i.e., the complex (\ref{v66}) is one-exact.
\end{proof}

\section{The Koszul--Tate complex}

Turn now to the second homology $H_2(\dl_0)$ of the complex
(\ref{v66}). A generic two-chain  reads
\mar{v77}\beq
\Phi= G + H= \op\sum_{0\leq|\La|} G_r^\La c^r_\La  +
\op\sum_{0\leq|\La|,|\Si|} H_{ab}^{\La\Si} c^a_\La c^b_\Si, \qquad
G_r^\La, H_{ab}^{\La\Si} \in \cO_\infty Y. \label{v77}
\eeq
The cycle condition $\dl_0 \Phi=0$ takes the form
\mar{v79}\beq
 \op\sum_{0\leq|\La|} G_r^\La d_\La\Delta^r +\dl H=0.
\label{v79}
\eeq
One can think of this equality as being the reduction condition on
the Noether operators $\Delta^r$ (\ref{v71}). Conversely, let
\be
\Phi=\op\sum_{0\leq|\La|} G_r^\La c^r_\La\in \cP_\infty[E;Y;E_0]_2
\ee
be a graded function such that the reduction condition (\ref{v79})
holds. Obviously, it is a cycle condition of the two-chain
(\ref{v77}). The reduction condition (\ref{v79}) is trivial either
if a two-cycle $\Phi$ (\ref{v77}) is a boundary or its summand $G$
vanishes on $\Ker\cE$.

\begin{defi} \label{v133} \mar{v133}
A degenerate differential operator is said to be one-stage
reducible if there exist non-trivial reduction conditions
(\ref{v79}), called first-stage Noether identities.
\end{defi}

\begin{prop} \label{v134} \mar{v134}
First-stage Noether identities can be identified to nontrivial
elements of the homology $H_2(\dl_0)$ iff any $\dl$-cycle $\Phi\in
\cS_\infty[E;Y]_2$ is a $\dl_0$-boundary.
\end{prop}

\begin{proof}
It suffices to show that, if the summand $G$ of a two-cycle $\Phi$
(\ref{v77}) is $\dl$-exact, then $\Phi$ is a boundary. If $G=\dl
\Psi$, then
\mar{v169}\beq
\Phi=\dl_0\Psi +(\dl-\dl_0)\Psi + H. \label{v169}
\eeq
The cycle condition reads
\be
\dl_0\Phi=\dl((\dl-\dl_0)\Psi + H)=0.
\ee
Then $(\dl-\dl_0)\Psi + H$ is $\dl_0$-exact since any $\dl$-cycle
$\f\in \cS_\infty[E;Y]_2$ by assumption is a $\dl_0$-boundary.
Consequently, $\Phi$ (\ref{v169}) is $\dl_0$-exact. Conversely,
let $\Phi\in \cS_\infty[E;Y]_2$ be an arbitrary $\dl$-cycle.
 The cycle condition reads
\mar{v100}\beq
\dl\Phi= 2\Phi_{ab}^{\La\Si}c^a_\La\dl c^b_\Si=
2\Phi_{ab}^{\La\Si} c^a_\La d_\Si \cE^b=0. \label{v100}
\eeq
It follows that $\Phi_{ab}^{\La\Si} \dl c^b_\Si=0$ for all induces
$(a,\La)$. Omitting a $\dl$-boundary term, we obtain
\be
\Phi_{ab}^{\La\Si} c^b_\Si= G_{ar}^{\La\Xi}d_\Xi \Delta^r.
\ee
Hence, $\Phi$ takes the form
\mar{v135}\beq
\Phi=-G^{\La\Xi}_{ar} d_\Xi\Delta^r c^a_\La. \label{v135}
\eeq
We can associate to it the three-chain $\Psi= -G^{\La\Xi}_{ar}
c^r_\Xi c^a_\La$ such that
\be
\dl_0\Psi=\Phi +\si = \Phi - G^{\La\Xi}_{ar}d_\La\cE^a c^r_\Xi.
\ee
Owing to the equality $\dl\Phi=0$, we have $\dl_0\si=0$. Since
$\si$ is $\dl$-exact, it by assumption is $\dl_0$-exact, i.e.,
$\si=\dl_0\psi$.  Then we obtain that $\Phi=\dl_0\Psi -\dl_0\psi$.
\end{proof}

\begin{lem} \label{v200} \mar{v200}
It is easily justified that a two-cycle $\Phi\in
\cS_\infty[E;Y]_2$ is  $\dl_0$-exact iff $\Phi$ up to a
$\dl$-boundary takes the form
\mar{v140}\beq
\Phi= \op\sum_{0\leq |\La|, |\Si|} G^{\Si\La}_{rr'} d_\Si\Delta^r
d_\La\Delta^{r'}. \label{v140}
\eeq
\end{lem}

If the condition of Proposition \ref{v134} (called the
two-homology regularity condition) is satisfied, let us assume
that first-stage Noether identities are finitely generated as
follows. There are elements $\Delta_{(1)}\in H_2(\dl_0)$ making up
a projective $C^\infty(X)$-module $\cC_{(1)}$ of finite rank which
is isomorphic to the module of sections of some vector bundle
$E_1\to X$. Let $\{\Delta^{r_1}\}$ be local bases for this
$C^\infty(X)$-module. Every element $\Phi\in H_2(\dl_0)$
factorizes
\mar{v80,1}\ben
&& \Phi= \op\sum_{0\leq|\Xi|} \Phi_{r_1}^\Xi d_\Xi \Delta^{r_1},
\qquad \Phi_{r_1}^\Xi\in
\cO_\infty Y, \label{v80}\\
&&\Delta^{r_1}=G^{r_1}+ h^{r_1}=\op\sum_{0\leq|\La|} \Delta^{\La
r_1}_r c^r_\La + h^{r_1}, \qquad
 h^{r_1}\in
\cS_\infty[E;Y]_2, \label{v81}
\een
via elements of $\cC_{(1)}$, i.e., any first-stage Noether
identity (\ref{v79}) results from the equalities
\mar{v82}\beq
 \op\sum_{0\leq|\La|} \Delta^{r_1\La}_r d_\La \Delta^r +\dl
h^{r_1} =0, \label{v82}
\eeq
called the complete first-stage Noether identities. Elements of
$\cC_{(1)}$ are said to be the first-stage Noether operators. Note
that the first summands $G^{r_1}$ of the operators $\Delta^{r_1}$
(\ref{v81}) are not $\dl$-exact.

\begin{prop} \label{v139} \mar{v139} Given a reducible degenerate
differential operator $\cE$, let the associated one-exact complex
(\ref{v66}) obey the two-homology regularity condition and let its
homology $H_2(\dl_0)$ be finitely generated. Then this complex is
extended to the two-exact one with a boundary operator whose
nilpotency conditions are equivalent to complete Noether and
first-stage Noether identities (see the complex (\ref{v87})
below).
\end{prop}

\begin{proof}
Let us consider the ring $\cP^*_\infty[E_1E;Y;E_0]$ of graded
functions on $J^\infty Y$ possessing local bases
$\{y^i_\La,c^a_\La,c^r_\La,c^{r_1}_\La\}$, where $[c^{r_1}_\La]=1$
and Ant$[c^{r_1}_\La]=3$. It can be provided the first-stage
Koszul--Tate differential defined as the nilpotent graded
derivation
\mar{v205}\beq
\dl_1=\dl_0 + \rdr_{r_1} \Delta^{r_1}. \label{v205}
\eeq
Its nilpotency conditions (\ref{0688}) are equivalent to complete
Noether identities (\ref{v64}) and complete first-stage Noether
identities (\ref{v82}). Then the module
$\cP_\infty[E_1E;Y;E_0]_{\leq 4}$ of graded functions of antifield
number Ant$[\f]\leq 4$ is split into the chain complex
\mar{v87}\ben
&& 0\lto \im\dl \llr^\dl \cS_\infty[E;Y]_1\llr^{\dl_0}
\cP_\infty[E;Y;E_0]_2\llr^{\dl_1} \cP_\infty[E_1E;Y;E_0]_3
\label{v87}\\
&& \qquad \llr^{\dl_1} \cP_\infty[E_1E;Y;E_0]_4.\nonumber
\een
Let $H_*(\dl_1)$ denote its homology. It is readily observed that
\be
H_0(\dl_1)=H_0(\dl)=0, \qquad H_1(\dl_1)=H_1(\dl_0)=0.
\ee
By virtue of the expression (\ref{v80}), any two-cycle of the
complex (\ref{v87}) is a boundary
\be
 \Phi= \op\sum_{0\leq|\Xi|} \Phi_{r_1}^\Xi d_\Xi \Delta^{r_1}
=\dl_1(\op\sum_{0\leq|\Xi|} \Phi_{r_1}^\Xi c_\Xi^{r_1}).
\ee
It follows that $H_2(\dl_1)=0$, i.e., the complex (\ref{v87}) is
two-exact.
\end{proof}

If the third homology $H_3(\dl_1)$ of the complex (\ref{v87}) is
not trivial, there are reduction conditions on the first-stage
Noether operators, and so on. Iterating the arguments, we come to
the following.

Let $\cE$ be a degenerate differential operator whose Noether
identities are finitely generated. In accordance with Proposition
\ref{v137}, we associates to it the one-exact chain complex
(\ref{v66}). Given an integer $N\geq 1$, let $E_1, \ldots, E_N$ be
some vector bundles over $X$. Let us consider the ring
\be
\cP_\infty\{N\}=\cP_\infty[E_{N-1}\cdots E_1E;Y;E_0\cdots E_N]
\ee
of graded functions on $J^\infty Y$ if $N$ is even or the ring
\be
\cP_\infty\{N\}=\cP_\infty[E_N\cdots E_1E;Y;E_0\cdots E_{N-1}]
\ee
if $N$ is odd. It possesses local bases $\{y^i_\La, c^a_\La,
c^r_\La, c^{r_1}_\La, \ldots, c^{r_N}_\La\}$ where
$[c_\La^{r_k}]=(k+1)\,$mod2 and Ant$[c_\La^{r_k}]=k+2$. Let
$k=-1,0$ further stand for $y^i$ and $c^r$, respectively. We
assume that:

(i) the ring $\cP_\infty\{N\}$ is provided with the nilpotent
graded derivation
\mar{v92,'}\ben
&&\dl_N=\dl_0 + \op\sum_{1\leq k\leq N}\rdr_{r_k} \Delta^{r_k},
\label{v92}\\
&& \Delta^{r_k}=G^{r_k} + h^{r_k}= \op\sum_{0\leq|\La|}
\Delta^{\La r_k}_{r_{k-1}} c_\La^{r_{k-1}} + \op\sum_{0\leq \Si,
0\leq\Xi}(h^{\Xi\Si r_k}_{a r_{k-2}} c^a_\Xi c_\Si^{r_{k-2}}+...),
\label{v92'}
\een
of antifield number -1;

(ii) the module $\cP_\infty\{N\}_{\leq N+3}$ of graded functions
of antifield number $k\leq N+3$ is split into the $(N+1)$-exact
chain complex
\mar{v94}\ben
&&0\lto \im \dl \llr^\dl \cS_\infty[E;Y]_1\llr^{\dl_0}
\cP_\infty\{0\}_2\llr^{\dl_1} \cP_\infty\{1\}_3\cdots
\label{v94}\\
&& \qquad
 \llr^{\dl_{N-1}} \cP_\infty\{N-1\}_{N+1}
\llr^{\dl_N} \cP_\infty\{N\}_{N+2}\llr^{\dl_N}
\cP_\infty\{N\}_{N+3}, \nonumber
\een
which satisfies the $(N+1)$-homology regularity condition,
introduced below.

\begin{defi} \label{v155} \mar{v155} One says that the complex (\ref{v94})
obeys the $(N+1)$-homology regularity condition if any
$\dl_{k<N-1}$-cycle $\Phi\in \cP_\infty\{k\}_{k+3}\subset
\cP_\infty\{k+1\}_{k+3}$ is a $\dl_{k+1}$-boundary.
\end{defi}

\begin{rem}
The $(N+1)$-exactness of the complex (\ref{v94}) implies that any
$\dl_{k<N-1}$-cycle $\Phi\in \cP_\infty^{0,n}\{k\}_{k+3}$, $k<N$,
is a $\dl_{k+2}$-boundary, but not necessary a $\dl_{k+1}$-one.
\end{rem}

If $N=1$, the complex $\cP_\infty\{1\}_{\leq 4}$ (\ref{v94})
restarts the complex (\ref{v87}) associated to a first-stage
reducible differential operator in accordance with Proposition
\ref{v139}. Therefore, we agree to call $\dl_N$ (\ref{v92}) the
$N$-stage Koszul--Tate differential. Its nilpotency implies
complete Noether identities (\ref{v64}), first-stage Noether
identities (\ref{v82}) and the equalities
\mar{v93}\beq
\op\sum_{0\leq|\La|} \Delta^{\La r_k}_{r_{k-1}}d_\La
(\op\sum_{0\leq|\Si|} \Delta^{\Si r_{k-1}}_{r_{k-2}}
c_\Si^{r_{k-2}}) + \dl(\op\sum_{0\leq \Si, 0\leq\Xi}h^{\Xi\Si
r_k}_{a r_{k-2}} c_\Xi^a c_\Si^{r_{k-2}})=0, \label{v93}
\eeq
for $k=2,\ldots,N$. One can think of the equalities (\ref{v93}) as
being complete $k$-stage Noether identities because of their
properties which we will justify in the case of $k=N+1$.
Accordingly, $\Delta^{r_k}$ (\ref{v92'}) are said to be the
$k$-stage Noether operators.

A generic $(N+2)$-chain $\Phi\in \cP_\infty\{N\}_{N+2}$ takes the
form
\mar{v156}\beq
\Phi= G + H= \op\sum_{0\leq|\La|} G_{r_N}^\La c_\La^{r_N} +
\op\sum_{0\leq \Si, 0\leq\Xi}(H^{\Xi\Si}_{a r_{N-1}} c_\Xi^a
c_\Si^{r_{N-1}}+...). \label{v156}
\eeq
Let it be a cycle. The cycle condition $\dl_N\Phi=0$ implies the
equality
\mar{v145}\beq
\op\sum_{0\leq|\La|} G_{r_N}^\La d_\La (\op\sum_{0\leq|\Si|}
\Delta^{\Si r_N}_{r_{N-1}}c_\Si^{r_{N-1}}) + \dl(\op\sum_{0\leq
\Si, 0\leq\Xi}H^{\Xi\Si}_{ar_{N-1}} c_\Xi^a c_\Si^{r_{N-1}})=0.
\label{v145}
\eeq
One can think of this equality as being the reduction condition on
the $N$-stage Noether operators (\ref{v92'}). Conversely, let
\be
\Phi= \op\sum_{0\leq|\La|} G_{r_N}^\La c_\La^{r_N} \in
\cP_\infty\{N\}_{N+2}
\ee
be a graded function such that the reduction condition
(\ref{v145}) holds. Then this reduction condition can be extended
to a cycle one as follows. It is brought into the form
\be
\dl_N(\op\sum_{0\leq|\La|} && G_{r_N}^\La c_\La^{r_N} +
\op\sum_{0\leq \Si, 0\leq\Xi}H^{\Xi\Si}_{ar_{N-1}}
c_\Xi^a c_\Si^{r_{N-1}})=\\
&& \qquad  -\op\sum_{0\leq|\La|} G_{r_N}^\La d_\La h^{r_N}
+\op\sum_{0\leq \Si, 0\leq\Xi}H^{\Xi\Si}_{ar_{N-1}} c_\Xi^a d_\Si
\Delta^{r_{N-1}}.
\ee
A glance at the expression (\ref{v92'}) shows that the term in the
right-hand side of this equality belongs to
$\cP_\infty\{N-2\}_{N+1}$. It is a $\dl_{N-2}$-cycle and,
consequently, a $\dl_{N-1}$-boundary $\dl_{N-1}\Psi$ in accordance
with the $(N+1)$-homology regularity condition. Then the reduction
condition (\ref{v145}) is a $c_\Si^{r_{N-1}}$-dependent part of
the cycle condition
\be
\dl_N(\op\sum_{0\leq|\La|} G_{r_N}^\La c_\La^{r_N} +
\op\sum_{0\leq \Si, 0\leq\Xi}H^{\Xi\Si}_{ar_{N-1}}c_\Xi^a
c_\Si^{r_{N-1}} -\Psi)=0,
\ee
but $\dl_N\Psi$ does not make a contribution to this reduction
condition.

Being a cycle condition, the reduction condition (\ref{v145}) is
trivial either if a cycle $\Phi$ (\ref{v156}) is a
$\dl_N$-boundary or its summand $G$ is $\dl$-exact, i.e., it is a
boundary, too, as we have stated above. Then Definition \ref{v133}
can be generalized as follows.

\begin{defi} \label{v202} \mar{v202}
A degenerate differential operator is said to be $(N+1)$-stage
reducible if there exist non-trivial reduction conditions
(\ref{v145}), called the $(N+1)$-stage Noether identities.
\end{defi}

\begin{theo} \label{v163} \mar{v163}
(i) The $(N+1)$-stage Noether identities can be identified to
nontrivial elements of the homology $H_{N+2}(\dl_N)$ of the
complex (\ref{v94}) iff this homology obeys the $(N+2)$-homology
regularity condition. (ii) If the homology $H_{N+2}(\dl_N)$ is
finitely generated as defined below, the complex (\ref{v94})
admits an $(N+2)$-exact extension.
\end{theo}

\begin{proof}
(i) The $(N+2)$-homology regularity condition implies that any
$\dl_{N-1}$-cycle $\Phi\in \cP_\infty\{N-1\}_{N+2}\subset
\cP_\infty\{N\}_{N+2}$ is a $\dl_N$-boundary. Therefore, if $\Phi$
(\ref{v156}) is a representative of a nontrivial element of
$H_{N+2}(\dl_N)$, its summand $G$ linear in $c_\La^{r_N}$ does not
vanish. Moreover, it is not a $\dl$-boundary. Indeed, if $\Phi=\dl
\Psi$, then
\mar{v172}\beq
\Phi=\dl_N\Psi +(\dl-\dl_N)\Psi + H. \label{v172}
\eeq
The cycle condition takes the form
\be
\dl_N\Phi=\dl_{N-1}((\dl-\dl_N)\Psi + H)=0.
\ee
Hence, $(\dl-\dl_0)\Psi + H$ is $\dl_N$-exact since any
$\dl_{N-1}$-cycle $\f\in \cP_\infty\{N-1\}_{N+2}$ is a
$\dl_N$-boundary.  Consequently, $\Phi$ (\ref{v172}) is a
boundary. If the $(N+2)$-homology regularity condition does not
hold, trivial reduction conditions (\ref{v145}) also come from
nontrivial elements of the homology $H_{N+2}(\dl_N)$.

(ii) Let the $(N+1)$-stage Noether identities be finitely
generated. Namely, there exist elements $\Delta_{(N+1)}\in
H_{N+2}(\dl_N)$ making up a projective $C^\infty(X)$-module
$\cC_{(N+1)}$ of finite rank which is isomorphic to the module of
sections of some vector bundle $E_{N+1}\to X$. Let
$\{\Delta^{r_{N+1}}\}$ be local bases for this
$C^\infty(X)$-module. Then any element $\Phi\in H_{N+2}(\dl_N)$
factorizes
\mar{v160,1}\ben
&& \Phi= \op\sum_{0\leq|\Xi|} \Phi_{r_{N+1}}^\Xi d_\Xi
\Delta^{r_{N+1}}, \qquad \Phi_{r_{N+1}}^\Xi\in
\cO_\infty Y, \label{v160}\\
&&\Delta^{r_{N+1}}=G^{r_{N+1}}+ h^{r_{N+1}}=\op\sum_{0\leq|\La|}
\Delta^{\La r_{N+1}}_{r_N} c_\La^{r_N} + h^{r_{N+1}}, \label{v161}
\een
through elements of $\cC_{(N+1)}$. Clearly, this factorization is
independent of specification of local bases
$\{\Delta^{r_{N+1}}\}$. Let us extend the ring $\cP_\infty\{N\}$
 to the  ring $\cP_\infty\{N+1\}$ possessing local bases
\be
\{y^i_\La, c_\La^a, c_\La^r, c_\La^{r_1}, \ldots, c_\La^{r_N},
c_\La^{r_{N+1}}\}, \qquad [c_\La^{r_{N+1}}]=N{\rm mod}\,2, \qquad
{\rm Ant}[c_\La^{r_{N+1}}]=N+3.
\ee
It is provided with the nilpotent graded derivation
\mar{v170}\beq
\dl_{N+1}=\dl_N + \rdr_{r_{N+1}} \Delta^{r_{N+1}} \label{v170}
\eeq
of antifield number -1. With this graded derivation, the module
$\cP_\infty\{N+1\}_{\leq N+4}$ of graded functions of antifield
number Ant$[f]\leq N+4$ is split into the chain complex
\mar{v171}\ben
&&0\lto \im \dl \llr^\dl \cS_\infty[E;Y]_1\llr^{\dl_0}
\cP_\infty\{0\}_2\llr^{\dl_1} \cP_\infty\{1\}_3\cdots
\llr^{\dl_{N-1}} \cP_\infty\{N-1\}_{N+1}
 \label{v171}\\
&& \qquad \llr^{\dl_N} \cP_\infty\{N\}_{N+2}\llr^{\dl_{N+1}}
\cP_\infty\{N+1\}_{N+3}\llr^{\dl_{N+1}} \cP_\infty\{N+1\}_{N+4}.
\nonumber
\een
It is readily observed that this complex is $(N+2)$-exact. In this
case, the $(N+1)$-stage Noether identities (\ref{v145}) come from
the complete $(N+1)$-stage Noether identities
\be
\op\sum_{0\leq|\La|} \Delta^{\La r_{N+1}}_{r_N} d_\La \Delta^{r_N}
+\dl h^{r_{N+1}} =0,
\ee
which are reproduced as the nilpotency conditions of the graded
derivation (\ref{v170}).
\end{proof}

The iteration procedure based on Theorem \ref{v163} can be
prolonged up to an integer $N_{\rm max}$ when the $N_{\rm
max}$-stage Noether identities are irreducible, i.e., the homology
$H_{N_{\rm max}+2}(\dl_{N_{\rm max}})$ is trivial. This iteration
procedure may also be infinite. It results in the manifested exact
Koszul--Tate complex with the Koszul--Tate boundary operator whose
nilpotency conditions reproduce all Noether and higher Noether
identities of an original differential operator $\cE$.


\begin{thebibliography}{ddd}

\bibitem{jpa05} D. Bashkirov, G. Giachetta, L. Mangiarotti,  and
G. Sardanashvily, Noether's second theorem in a general setting.
Reducible gauge theories, {\it J. Phys. A} {\bf 38} (2005),
5239-5344.

\bibitem{jmp05} D. Bashkirov, G. Giachetta, L. Mangiarotti, and
G. Sardanashvily, Noether's second theorem for BRST symmetries,
{\it J. Math. Phys.} {\bf 46} (2005), 053517.

\bibitem{kras} I.Krasil'shchik, V.Lychagin and A.Vinogradov, {\it Geometry of
Jet Spaces and Nonlinear Partial Differential Equations} (Gordon
and Breach, Glasgow, 1985).

\bibitem{bry} R.Bryant, S.Chern, R.Gardner, H.Goldschmidt and P.Griffiths,
{\it Exterior Differential Systems} (Springer-Verlag, Berlin,
1991).

\bibitem{cmp04} G. Giachetta, L. Mangiarotti, and G. Sardanashvily,
Lagrangian supersymmetries depending on derivatives. Global
analysis and cohomology, {\it Commun. Math. Phys.} (2005)
(accepted).

\bibitem{bart} C. Bartocci, U. Bruzzo, and D. Hern\'andez
Ruip\'erez, {\it The Geometry of Supermanifolds} (Kluwer,
Dordrecht, 1991).

\bibitem{book05} G. Giachetta, L. Mangiarotti, and G. Sardanashvily,
{\it Geometric and Algebraic Topological Methods in Quantum
Mechanics} (World Scientific, Singapore, 2005).

\bibitem{ren} A. Rennie, Smoothnes and locality for nonunital spectral triples,
 {\it K-Theory} {\bf 28} (2003), 127-165.

\bibitem{hern} D. Hern\'andez Ruip\'erez and J. Mu\~noz
Masqu\'e, Global variational calculus on graded manifolds, {\it J.
Math. Pures Appl.} {\bf 63} (1984), 283-309.

\bibitem{barn} G. Barnich, F. Brandt and M. Henneaux,  Local
BRST cohomology in gauge theories, {\it Phys. Rep.} {\bf 338}
(2000), 439-569.



\end{thebibliography}
\end{document}